\documentclass[10pt]{article}
\usepackage{amsmath, amsfonts, amsthm, amssymb, verbatim,dsfont,a4wide}
\usepackage{graphicx}
\usepackage[mathcal]{euscript}
\usepackage[applemac]{inputenc}
\usepackage{dsfont} 
\usepackage{amssymb,mathrsfs}
\usepackage[usenames]{color}

\usepackage[T1]{fontenc}%

\theoremstyle{plain}

\theoremstyle{definition}

\theoremstyle{plain}
        
\numberwithin{equation}{section}

\usepackage{amsmath}
\usepackage{verbatim}
\newcommand \be           {\begin{equation}}
\newcommand \ee            {\end{equation}}

\newcommand \RR           {\mathbb{R}}

\newcommand \Pbold           {\mathbf{P}} 
\newcommand \PP \Pbold

\newcommand \del           \partial
\newcommand \eps            \epsilon

\def\XXint#1#2#3{{\setbox0=\hbox{$#1{#2#3}{\int}$}
\vcenter{\hbox{$#2#3$}}\kern-.5\wd0}}

\let\oldmarginpar\marginpar
\renewcommand\marginpar[1]{\-\oldmarginpar[\raggedleft\footnotesize #1]%
{\raggedright\footnotesize #1}}
\def\build#1_#2^#3{\mathrel{
\mathop{\kern 0pt#1}\limits_{#2}^{#3}}}
\allowdisplaybreaks[1]

\begin{document}

\title{A continuous model of ant foraging\\ with pheromones and trail formation}
\author{Paulo Amorim$^1$}
\footnotetext[1]{
Instituto de Matem\'atica, Universidade Federal do Rio de Janeiro,
Av. Athos da Silveira Ramos 149,
Centro de Tecnologia - Bloco C,
Cidade Universit\'aria - Ilha do Fund\~ ao,
Caixa Postal 68530, 21941-909 Rio de Janeiro,
RJ - Brasil\\
Email: paulo@im.ufrj.br. web page: \texttt{http://www.im.ufrj.br/$\sim$paulo/}}
\date{}

\maketitle

\begin{abstract} 
We propose and numerically analyze a PDE model of ant foraging behavior. Ant foraging is a prime example of individuals following simple behavioral rules based on local information producing complex, organized and ``intelligent'' strategies at the population level. One of its main aspects is the widespread use of pheromones, which are chemical compounds laid by the ants used to attract other ants to a food source. In this work, we consider a continuous description of a population of ants and simulate numerically the foraging behavior using a system of PDEs of chemotaxis type. We show that, numerically, this system accurately reproduces observed foraging behavior, such as trail formation and efficient removal of food sources. 
\end{abstract}

\section{Introduction}

Ant foraging is among the most interesting emergent behaviors in the social insects, and indeed in the animal kingdom. Perhaps the most striking aspect of ant foraging is how individuals following simple behavioral rules based on local information
produce complex, organized and seemingly intelligent strategies at the population level. Ant foraging (along with most other activities of an ant colony) is a prime example of so-called \emph{emergent} behavior.  Indeed, while
individuals of an ant population possess only very limited cognitive abilities, each is equipped with a basic set of behavioral rules, to the effect that the emerging collective behavior exhibits a remarkable degree of efficiency, optimality, and adaptivity to a changing environment. 

It has long been known that one of the main forms of communication among ants is the use of pheromones. These are chemical compounds which 
individual ants secrete and leave on the substrate and which in effect are used as a means of communication between ants, transmitting a variety of messages such as alarm, presence of food, or providing colony-specific olfactory signatures used by the ants to identify nest-mates.

Among the many documented functions of ant 
pheromones, we are interested in their role as a chemical trail indicating the 
direction to a food source. Many species of ants, especially trail-forming 
ones, are known to lay a pheromone as they travel from the food source back to 
nest. The main attribute of this pheromone is that it is attractive to other 
ants, who tend to follow the direction of increasing concentration of the 
chemical. These ants will then reach the food source and return to the nest 
while laying pheromone themselves, thus reinforcing the chemical trail in a 
positive feedback loop. This results in the formation of well defined trails 
leading from the nest to the food source, allowing for an efficient transport 
of the food to the nest. 
Thus, pheromones play a major role in food foraging, where they are widely used (among other strategies) to recruit nest mates to new food sources.  

The entomological research body on ants, their behavior, and their olfactory means of communication is vast. Here, we content ourselves with citing some seminal works, as well as some more recent investigations with special relevance to our analysis. For a general reference on myrmecology (the branch of entomology that deals with ants), we refer to the encyclopedic book by H\"olldobler and Wilson \cite{TheAnts}. Therein may be found many relevant references up to 1990. The paper \cite{Regnier1968} contains an overview of the chemical study of pheromones. For some modern references on the trail-laying behavior of ants, we refer to \cite{Bossert1963,Deneubourg1990,Beckers1992,Edelstein1995,Sumpter2003,Vittori2004,Van1986}, and the references therein, although of course many other papers could be cited.
Concerning the computational simulation of ant trail-laying, we refer to \cite{Boissard2012,Couzin2002,Edelstein1995,Johnson2006,Schweitzer1997,Sumpter2003,Watmough1995,Watmough1995-2,Weyer1985}, though again this list is far from complete. See especially \cite{Boissard2012} for a recent PDE approach (though quite different from our own), and an in-depth literature review of available numerical and modeling strategies for ant foraging.

In particular, from what our bibliographical research could gather, only the work \cite{Watmough1995} presents a PDE model which (as our own) divides the ant population into two kinds, namely ants leaving the nest and ants returning to the nest. However, the setting in \cite{Watmough1995} is highly simplified, being one-dimensional, so no trail formation occurs. Moreover, the system proposed in \cite{Watmough1995} is only explored numerically in a simplified ODE version.

Thus, to the best of our knowledge, the present work is the first to consider the modeling and simulation of the whole cycle of food foraging by ants, comprising random foraging, discovery and transport of food, formation of trails and fading of trails upon exhaustion of the food sources.

\section{Modeling foraging in ant colonies: a continuous model}
In this Section we present our PDE model of ant foraging behavior. Our model incorporates successfully the following traits:
\begin{itemize}
\item spontaneous trail formation in the presence of food sources,
\item discovery of new food sources,
\item depletion of food sources,
\item abandoning of unproductive trails, and
\item topography.
\end{itemize}

\subsection{Modeling assumptions}
We intend to model not one specific species of ant, but rather to capture in a qualitative way the characteristic properties of ant foraging. In this sense, one assumption we make in our model is that the ants know the way back to their nest upon finding a food source. This is reflected by the introduction of a given nest-bound vector field $\mathbf v(x)$. Importantly, this assumption is supported by the literature. Indeed, some species of ants have been proven to use visual and olfactory orientation cues to return to the nest \cite{Holldobler1976,Steck2009}, as well as so-called orientation by path integration \cite{Muller1988}, in which individuals cumulatively keep track of changes in direction and thus of the overall direction of the nest.
Therefore the introduction of the nest-bound field $\mathbf v(x)$ is realistic but still simplifies our algorithms and equations.

Of course, not all species of ants will employ one of the above methods for returning to the nest. In the myrmecological literature, it is difficult to find discussions of the means by which ants orient themselves towards the nest. This is in contrast to the well studied methods that ants use to find food.

\subsection{Continuous description of foraging behavior.}
In this work, we study ant foraging behavior from the mathematical point of view of \emph{chemotaxis} (see \cite{HillenPainter2009} for an up-to-date review). The term chemotaxis is used to describe any phenomena in which the movement of an agent (usually a cell or bacteria) is affected by the presence of a chemical agent. Typically, individuals follow paths of increasing concentration of the chemical agent, and often produce the agent themselves thus originating a variety of phenomena, including finite time blow-up, segregation of species, and the formation of patterns \cite{HillenPainter2009}.

It is therefore natural (as had already been observed in \cite{Boissard2012}) to approach ant foraging behavior from a chemotactical point of view. Thus, the ant population is modeled by a density function rather than by a discrete set of individuals. The ants are dynamically divided into two groups, as follows:
\begin{itemize}
\item Individuals leaving the nest randomly forage for food until they encounter food. Ants not carrying food (i.e., foraging) are called \emph{foraging ants}.
\item In the foraging phase, ants follow a random path (modeled by a diffusion term) and follow the gradient of the pheromone (if any). This predisposition to follow the pheromone is modeled by a transport term. Thus, when encountering a trail, they will tend to move in the direction of increasing density of pheromone, all the while avoiding steep topographical gradients.
\item When the ants encounter food, they become \emph{returning ants} which carry food back to the nest. We assume that ants know their way back to the nest, as explained above. As the ants return to the nest, they continuously lay  pheromone, which attracts foraging ants, but are not themselves attracted to it. 
\item Upon reaching the nest, returning ants transform back to foraging ants and leave the nest to continue foraging.
\item The pheromone evaporates at an exponential rate and diffuses according to a standard diffusion equation, while the food is depleted according to the quantity of foraging ants reaching it.
\end{itemize}

The population of foraging and returning ants is modeled by a \emph{density functions} $u(x,t)$ and $w(x,t)$ depending on $x\in \Omega$ and $t\ge 0$, where $\Omega$ is a bounded open subset of $\RR^2$, representing the physical domain that the population inhabits, and $t$ is time. Even though individual ants are not microscopic, we argue that modeling an ant population using a continuous density is reasonable. Indeed, to take the example of the genus \emph{Pogonomyrmex,} a typical worker measures 1.8mm in body length, while their foraging range usually extends to distances of 45--60m \cite{Holldobler1976}. Thus, the ratio of foraging distance to average body length is on the order of $3\times 10^4.$ The same assumption is used frequently, for instance, in the continuous modeling of cell dynamics \cite{HillenPainter2009}. In such settings, the typical ratio of cell size to physical domain is in the same order of magnitude as in the framework proposed here.

Also, we will see in \cite{Amorim} how the model below may be rigorously justified from a simple kinetic model which is in accordance with our modeling assumptions, as is common practice in the chemotaxis literature \cite{HillenPainter2009}.

Moreover, as is customary when modeling physical phenomena using reaction--diffusion equations, such as crowd movements or chemotaxis, the solution may be seen as the averaged outcome of a great number of individual experimental runs. 

\subsection{A simplified model.}
We begin by defining the variables and quantities used in our modeling framework:
\begin{itemize}
\item $u(x,t)=$ density of foraging ants
\item $w(x,t)=$ density of ants returning to the nest with food
\item $c(x,t)=$ concentration of food source
\item $v(x,t)=$ concentration of food pheromone
\item $z(x)=$ topography of terrain
\item $\mathbf v(x)=$ nest-bound field.
\end{itemize}

We propose a continuous model which can be described in words as:

\begin{itemize}
\item the foraging $u$-ants emerge from the nest located at $x=0$ at rate $M$ until time $T$ (or else their initial distribution may be a given function). They disperse according to Fourier's heat law ($\Delta u$ term), 
avoiding steep topography. Upon encountering the pheromone $v$, they follow its gradient ($\nabla\cdot (u\nabla v)$ term). 
\item When they reach an area with food (described by the food concentration 
$c$), they  turn into $w$-ants (returning ants), by means of the coupling source terms, depleting the food in the process. 
\item the $w$-ants now return to the nest with food, by following the (prescribed) home-bound vector field $\mathbf v$. They lay the pheromone $v$ which the $u$-ants will now follow to reach the food, according to the chemotactical transport term $\nabla\cdot( u\, \nabla v)$.
\item when the $w$-ants reach the nest located at $x=0$, they leave the food at the nest and re-emerge as $u$-ants. 
\item the third equation represents the laying of the pheromones by the $w$-ants, which evaporates and diffuses. 
\item the last equation represents the depletion of the food.
\end{itemize}

We now present a mathematical realization of the preceding description. This first simplified version contains every relevant structural aspect, which we will later refine through a finer modeling. The system reads, in nondimensional form (the details may be found in \cite{Amorim}),

\be
\label{1000}
\left\{
\begin{aligned}
&\partial_t u -  \nabla\cdot\big(\alpha_1 {\nabla u} + \alpha_2 u \, \nabla v  \big) = F(u,c,w)
\\
& \partial_t w - \nabla \cdot \big ( \alpha_3 \nabla w -  \alpha_2 \mathbf{v} w  \big) = G(u,c,w)
\\
& \partial_t v =  w -  v +  \Delta v
\\
& \partial_t c = - \alpha_4 u\, c
\\
&+ \text{Appropriate boundary conditions and initial data,}
\end{aligned} \right.
\ee
where the source terms $F$ and $G$ are found in \eqref{2000} below. These terms model ants exiting and entering the nest, as well as the transformation of foraging ants into returning ants.

\subsection{A Complete model}
Here we propose a more complete system for the description of ant foraging with pheromones and trail formation. It is based on the simplified system \eqref{1000} but incorporates some additional assumptions designed to render \eqref{1000} more realistic. So, in the presence of pheromone, the foraging ants should disperse less, whence a term $e^{-\alpha_3 v}$ limiting the diffusion of $u$. Also, a topography $z(x)$ is introduced in the transport term. Both foraging and returning ants try to avoid steep topography, which is modeled by the introduction of the terms $\nabla z$. Additionally, the terms $\beta(u),\beta(w)$ act to prevent overcrowding by limiting the density of ants, see \cite{HillenPainter2001}. The complete system in nondimensional form reads
\be
\label{2000}
\left\{
\begin{aligned}
&\partial_t u -  \nabla\cdot\big( \alpha_1 e^{-\gamma v}\nabla u + (   \nabla z - \alpha_2\nabla v  )\, u\beta(u) \big) = - uc + \alpha_5 w \delta_{x=0} + M(t) \delta_{x=0}
\\
& \partial_t w - \nabla \cdot \big ( \alpha_3 \nabla w + (  \nabla z - \alpha_2 \mathbf{v}) w \beta(w) \big) =  uc - \alpha_5 w \delta_{x=0}
\\
& \partial_t v =  w -  v +  \Delta v
\\
& \partial_t c = - \alpha_4 u\, c
\\
&+ \text{Appropriate boundary conditions and initial data.}
\end{aligned} \right.
\ee
Note that the source terms need not contain Diracs. Indeed, for the numerical simulations, we have used a small but extended nest, which is actually more realistic.

\section{Numerical results}
We present some results of the numerical simulation of system \eqref{2000}. We have plotted the density of the returning ants $w(x,t)$ in \eqref{2000} at successive times. The nest is at the center of each plot, and two food sources are placed at the ends of the observed trails, one above and one below the nest. One can clearly see the formation of trails emanating from the food sources to the nest, and that once the upper food source is depleted (at some moment between the third and fourth picture), the trail leading to that food source fades away due to the evaporation of the pheromone.

We have used a first order finite volume method on a square grid with standard Euler time-stepping. For the simulations presented here, we used a time step of 0.0005 for 300,000 iterations with 40,000 spatial cells.

$$\includegraphics[width=0.5\linewidth,keepaspectratio=false]{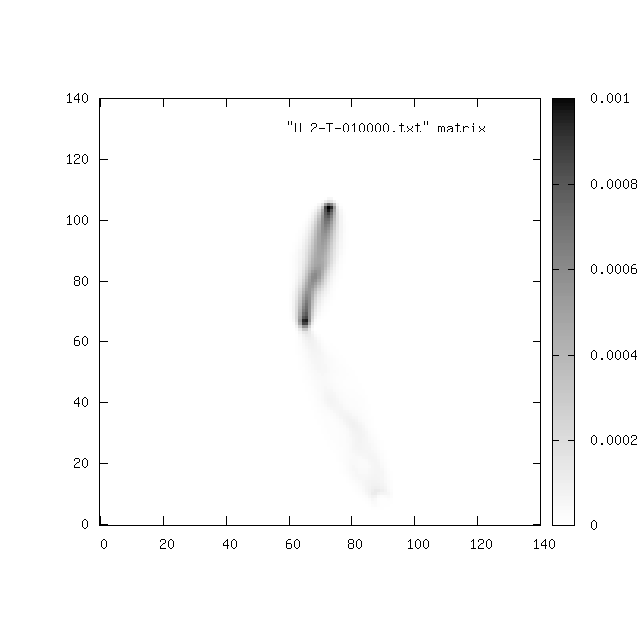}
\includegraphics[width=0.5\linewidth,keepaspectratio=false]{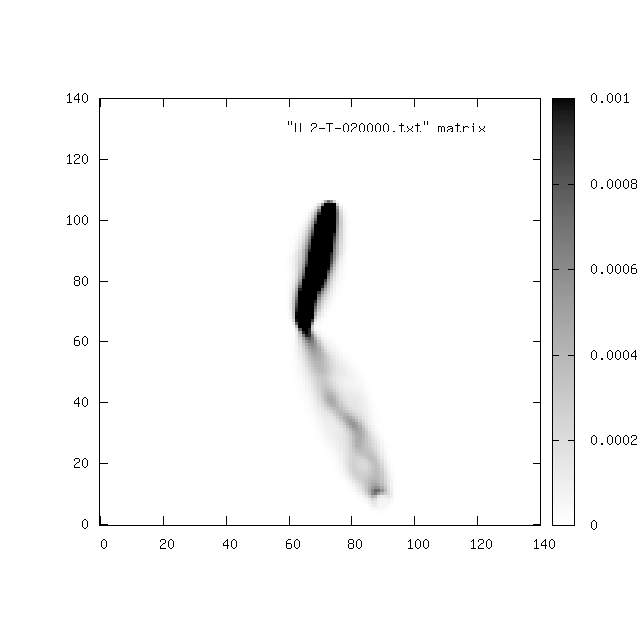}
$$
$$
\includegraphics[width=0.5\linewidth,keepaspectratio=false]{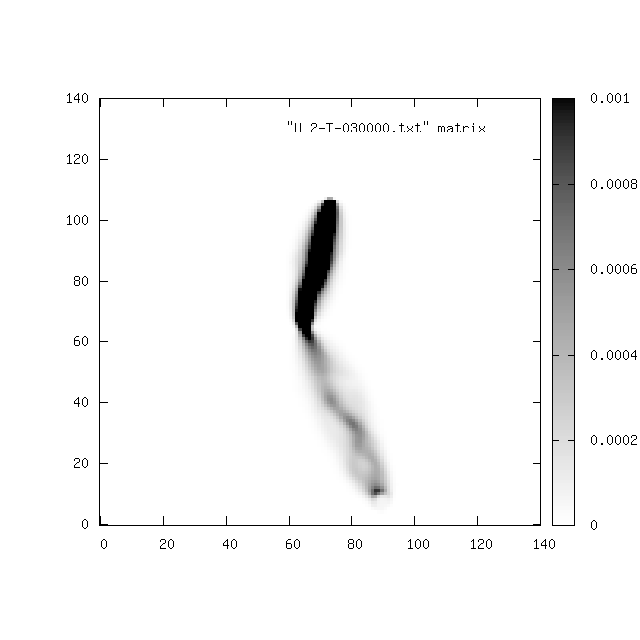}
\includegraphics[width=0.5\linewidth,keepaspectratio=false]{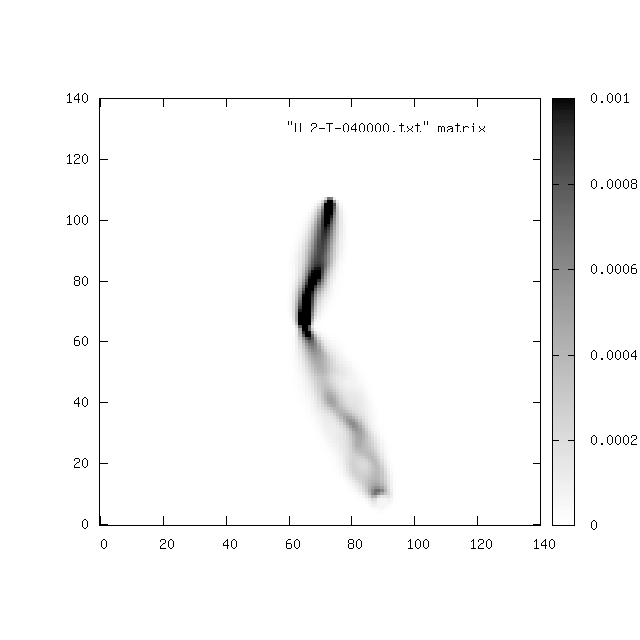}
$$
$$
\includegraphics[width=0.5\linewidth,keepaspectratio=false]{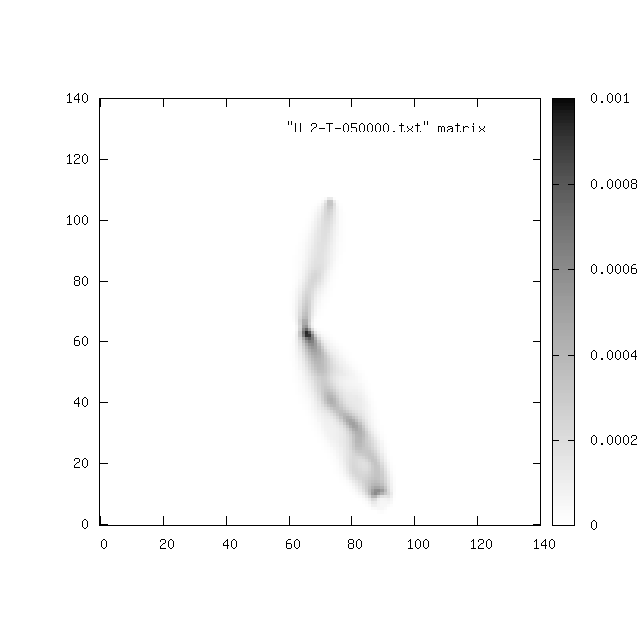}
\includegraphics[width=0.5\linewidth,keepaspectratio=false]{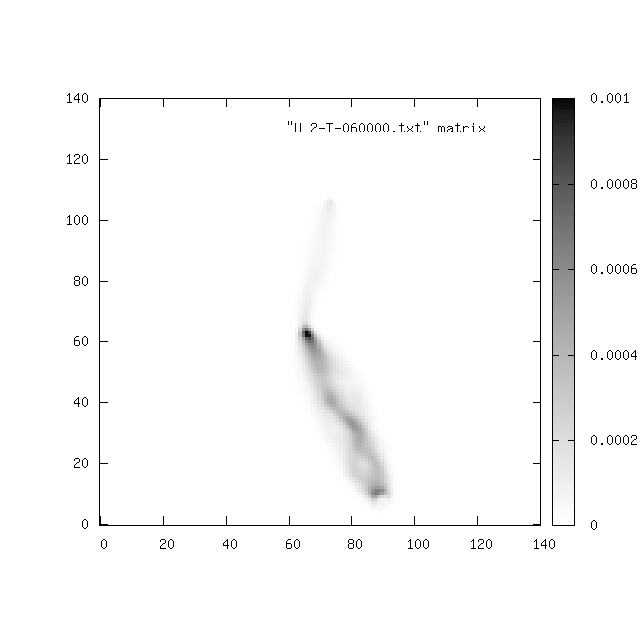}
$$

\section{Conclusions}
We have presented a mathematical model of ant foraging behavior using a system of PDEs in the framework of chemotaxis. We have shown numerically that this system exhibits spontaneous trail formation in the presence of food sources. Our model allows for the simulation of a whole cycle of food foraging, from ants emerging from a nest onto a foraging ground where food is placed, discovering the food and returning to the nest laying pheromones. Recruitment then takes place, with foraging ants being attracted in the vicinity of the nest to the pheromones laid previously by the returning ants. A feedback loop ensues, as more ants reach the food source and return to the nest laying pheromones. Finally, when the food source is exhausted, the trail fades away due to the natural evaporation of the pheromone.


%

\end{document}